\documentclass[11pt]{amsart} 
\pdfoutput=1
\usepackage{amssymb,amsmath,amsthm,amsfonts,amsbsy,latexsym,dsfont}

\newtheorem{prop}{Proposition}
\newtheorem{lem}{Lemma}
\newtheorem{thm}{Theorem}
\newtheorem{cor}{Corollary}

\newcommand{\be}{\begin{equation}}
\newcommand{\ee}{\end{equation}}

\newcommand{\real}{\mathbb{R}}

\newcommand{\N}{\mathbb{N}}

\newcommand{\Z}{\mathbb{Z}}

\newcommand{\C}{\mathcal{C}}

\newcommand{\calS}{\mathcal{S}}

\newcommand{\X}{\mathcal{X}}

\begin{document}

\bibliographystyle{plain}

\title[Uniform L-1 Ergodic Theorem]
{Entropy and the Uniform Mean Ergodic Theorem for a Family of Sets}

\date{March 2014}

\author[Adams]{Terrence M.\ Adams}
\address{Terrence Adams is with the Department of Defense, 
9161 Sterling Drive, Laurel, MD 20723} 
\author[Nobel]{Andrew B.\ Nobel} 
\address{Andrew Nobel is with the Department of Statistics and 
Operations Research, University of North Carolina at Chapel Hill, Chapel Hill, 
NC 27599-3260}
\thanks{The work of Andrew Nobel was supported by NSF Grants 
DMS-0907177 and DMS-1310002}
\email{nobel@email.unc.edu}

\begin{abstract} 
We define a notion of entropy for an infinite family $\C$ of measurable sets in a probability space.
We show that the mean ergodic theorem holds uniformly for $\C$ under every ergodic
transformation if and only if $\C$ has zero entropy.  
When the entropy of $\C$ is positive, we establish a strong converse showing that the uniform mean ergodic
theorem fails generically in every isomorphism class, including the isomorphism classes of 
Bernoulli transformations.   As a corollary of these results, we establish that every strong 
mixing transformation is uniformly strong mixing on $\C$ if and only if the entropy of $\C$ is zero, and  
obtain a corresponding result for weak mixing transformations.
\end{abstract}

\maketitle

\section{Introduction}

\noindent 
Let $T$ be an ergodic measure preserving transformation of a Lebesgue probability space 
$(\X, \calS, \mu)$.  For any fixed measurable set $C$, the 
ergodic theorem guarantees the pointwise and mean convergence of the relative
frequency $n^{-1} \sum_{i=0}^{n-1} I_C(T^i x)$ to $\mu(C)$.   
In many cases of theoretical and practical interest, it is useful (or sometimes essential) to know that 
the ergodic theorem holds uniformly over an infinite family $\C$ of measurable sets.  
A uniform pointwise
ergodic theorem holds for a family $\C$ if
\be
\label{upet}
\sup_{C \in \C} \, \Big| \frac{1}{n} \sum_{i=0}^{n-1} I_C(T^i x) - \mu(C) \Big| 
\ \to \ 
0 \ \mbox{ with probability one as $n$ tends to infinity.}
\ee
In probability and statistics, the trajectory $\{ T^{i} x : i \geq 0 \}$ in (\ref{upet}) is replaced by a 
stationary sequence $\{ X_i : i \geq 0 \}$ of random variables; trajectories constitute the
special case in which $X_i = T^i X_0$ where $X_0$ has distribution $\mu$.
Uniform pointwise ergodic theorems, also known as Glivenko-Cantelli
theorems, have been widely studied in the machine learning and empirical process
literature, most often in the case of independent identically distributed (i.i.d.) sequences of 
random variables; results and further references can be found in \cite{VaaWel96, Dudley99, Vapnik00}.  
There is also a companion literature concerning uniform 
pointwise ergodic theorems for dependent sequences of random variables, including the case of trajectories.
Adams and Nobel \cite{AdNob10a} established that (\ref{upet}) holds for every ergodic measure preserving 
transformation $T$ of a probability space
$(\X, \calS, \mu)$ provided that the family $\C$ has finite Vapnik-Chervonenkis (VC) dimension. 
(The definition of VC dimension is given in Section \ref{EZEF} below.)  Extensions 
and related results can be found in von Handel \cite{van13}; see \cite{AdNob10a,AdNob12,van13} and 
the references therein for more details.  

The focus of this paper is uniform mean ergodic theorems.  
A uniform mean ergodic theorem holds for a family of measurable sets $\C$ if
\be
\label{umet}
\sup_{C \in \C} \, \int \Big| \frac{1}{n} \sum_{i=0}^{n-1} I_C(T^i x) - \mu(C) \Big| \, d\mu 
\ \to \ 
0 \ \mbox{ as $n$ tends to infinity.}
\ee
Note that (\ref{upet}) clearly implies (\ref{umet}) but not vice versa.  For any i.i.d.\ sequence
of random variables, elementary arguments show that the uniform mean ergodic theorem 
holds for {\em any} family $\C$ of measurable sets, and the same holds true under the type
mixing conditions typically imposed in statistical problems.
Thus uniform mean ergodic theorems are primarily of interest for trajectories and related strongly 
dependent processes.  

We establish an equivalence between uniform mean ergodic theorems for
measure preserving transformations of $(\X, \calS, \mu)$ and the complexity of the family $\C$.
The complexity of $\C$ is quantified in terms of an entropy measure that we now define.
Given measurable sets $C_1, \ldots, C_n \in \calS$, let $\vee_{i=1}^{n} C_i$ denote their join,
or equivalently, the partition $\Pi$ of $\X$ having cells of the form 
$\tilde{C_1} \cap \cdots \cap \tilde{C_n}$ where each $\tilde{C}_i \in \{ C_i, C_i^c \}$.
The join $\Pi_0 \vee \Pi_1$ of two partitions is defined to be the join of their constituent
cells; thus $\Pi_0 \vee \Pi_1$ refines $\Pi_0$ and $\Pi_1$.
Recall that the entropy of finite measurable partition $\Pi$ of $\X$ is defined by
\[
H_\mu (\Pi) 
\, = \, 
-\sum_{A \in \Pi} \mu (A) \log{\mu (A)} .
\]
Here and in what follows we assume, without loss of generality, that all logarithms are base 2, 
and that each cell of a partition has positive measure.
Let $\C \subseteq \calS$ be a family of measurable subsets of $\X$.  For every
ordered sequence $\mathbb{S}=( C_1, C_2, \ldots )$ with $C_i \in \C$ let 
\be
\label{entseq}
H_\mu(\mathbb{S}) \, = \, \liminf_{n \to \infty} \frac{1}{n} H_\mu(\vee_{i=1}^{n} C_i) .
\ee
The {\em entropy} of the family $\C$ is defined by
\[
H_\mu(\C) = \mbox{sup} \, H_\mu(\mathbb{S}) ,
\]
where the supremum is over all ordered sequences $\mathbb{S}$ of sets in $\C$.
Our principal result is the following theorem.

\begin{thm}
\label{mainthm} 
Let $(\X, \calS, \mu)$ be a Lebesgue probability space,
and let $\C \subseteq \calS$ be a family of measurable subsets of $\X$.  
Then $H_\mu(\C)=0$ if and only if (\ref{umet}) holds
for every ergodic measure preserving 
transformation $T$~of $(\X,\calS,\mu)$. 
\end{thm}

The only-if part of Theorem \ref{mainthm} follows from elementary approximation properties of zero entropy 
families, and is valid for any probability space.   The converse part of the theorem requires us to 
exhibit a $\mu$-preserving ergodic tranformation $T$ for which (\ref{umet}) fails when $H_\mu(\C)$ is positive.
A stronger, generic, version of the converse is described below. 

\vskip.1in

\noindent
{\bf Remark:} Theorem \ref{mainthm} 
may also be expressed in the terminology
of stochastic processes.  Let $(\X, \calS, \mu)$ be a Lebesgue probability space
and let $\C \subseteq \calS$ be a family of measurable subsets of $\X$.  
It follows from Theorem \ref{mainthm} that $H_\mu(\C)=0$ if and only if 
\[
\sup_{C \in \C} \, \Big| \frac{1}{n} \sum_{i=0}^{n-1} I_C( U_i ) - \mu(C) \Big| 
\ \to \ 
0 
\ \mbox{ with probability one as $n \to \infty$}
\]
for every stationary ergodic process $U_0, U_1, \ldots$ taking values in 
$(\X, \calS)$ and such that $U_i$ has distribution $\mu$.

\vskip.2in

\subsection{Examples of Zero Entropy Families}  
\label{EZEF}
$  $

\vskip.1in

\noindent
{\bf Finite dimensional families.}
The (dual) Vapnik-Chervonenkis (VC) dimension of a family $\C$ is the largest $k \geq 1$
for which there exists sets $C_1,\ldots,C_k \in \C$ whose join $\vee_{i=1}^k C_i$ has 
cardinality $2^k$, the largest possible value for $k$ sets.  If $\C$ contains arbitrary large finite collections
with maximal joins, then its dimension is infinite.  It is known (cf.\ \cite{VaaWel96}) that if $\C$ has finite VC-dimension
then for $C_1, \ldots, C_n \in \C$ the
cardinality of $\vee_{i=1}^n C_i$ is bounded by a polynomial in $n$ that depends only on the dimension.  
It then follows from the elementary bound
$H_\mu(\Pi) \leq \log |\Pi|$ that any family $\C$ with finite VC-dimension
has entropy zero {\it for all} measures $\mu$.  Although it is considerably more indirect, one may reach 
the same conclusion by combining Theorem \ref{mainthm} with the principal result of \cite{AdNob10a}.
For general $\X$, the family of sets $\{ x : g(x) \geq 0 \}$ where $g$
ranges over a finite dimensional vector space of real valued functions on $\X$ has finite VC dimension. 
For $\X = \real^d$, families with finite VC dimension include convex polytopes with at most $k$ faces and
the family of open balls (with arbitrary center and radius).  
See \cite{VaaWel96} and \cite{Dudley99} for more details.
Note that zero entropy families may have infinite
dimension: it is easy to construct $\C = \{C_1, C_2, \ldots \}$
with infinite VC dimension such that $\mu(C_i) \to 0$.  

\vskip.1in

\noindent
{\bf Bracketing families.}  A family $\C$ is said to be bracketing with respect to a measure $\mu$ if for every $\delta > 0$
there exists a finite family ${\mathcal D}$ such that for each $C \in \C$ there exist sets
$U, V \in {\mathcal D}$ such that $U \subseteq C \subseteq V$
and $\mu(V \setminus U) < \delta$.  If $\C$ is bracketing with respect to $\mu$ one may readily show using the
characterization of zero-entropy families in Lemma \ref{entlem3} below that $H_\mu(\C) = 0$.   Alternatively, if
$\C$ is bracketing with respect to $\mu$ then it is easy to show that the uniform pointwise ergodic theorem
(\ref{upet}) holds for every $\mu$-preserving transformation of $(\X, \calS)$ and therefore $H_\mu(\C) = 0$ by 
Theorem \ref{mainthm}.  If $\X = \real^d$ and $\mu$ is absolutely continuous with respect to Lebesgue measure then
the family of convex subsets of $\X$ is bracketing with respect to $\mu$.  

\vskip.1in

\noindent
{\bf Clustered families.} It is easy to see that a family $\C$ has entropy $H_\mu(\C) = 0$ if there is a 
fixed measurable set $C_0$ such that
for each $\delta > 0$ at
most finitely many sets in $\C$ satisfy $\mu(C \triangle C_0 ) > \delta$.  

\vskip.1in

\noindent
{\bf Stability Properties.} Zero entropy families are closed under natural set theoretic operations.   A finite union of zero
entropy families has zero entropy.  If families $\C_1$ and $\C_2$ have zero entropy, then so too do the families of
complements $\{ C^c : C \in \C_1 \}$, unions $\{ C_1 \cup C_2 : C_1 \in \C_1, C_2 \in \C_2 \}$, and intersections
$\{ C_1 \cap C_2 : C_1 \in \C_1, C_2 \in \C_2 \}$.

\vskip.2in

\noindent
{\bf Remark:}
Recall that the entropy of a fixed measure preserving transformation $T$ is defined by
\be
\label{enttrans}
H_{\mu}(T) \ := \ \sup_{\Pi} \lim_{n \to \infty} \, n^{-1} H( \vee_{i=0}^{n-1} \, T^{-i} \Pi),
\ee
where the supremum is taken over all finite measurable partitions $\Pi$ of $\X$.  
We note that the definitions of $H_\mu(\C)$ and $H_\mu(T)$ do not coincide: in (\ref{entseq}) we scale by
the number of sets, equivalently two-way partitions, that give rise to the join, while in (\ref{enttrans})
one scales according to the number of iterates of $T^{-1}$, without regard to the cardinality of
the partition $\pi$.

\subsection{Strong Converse and Corollaries}

Let $\Phi$ be the space of invertible measure preserving transformations 
of $(\X, \calS, \mu)$ endowed with the weak topology, namely $\phi_r \to \phi$ if, 
for all measurable sets $A \in \calS$,
\[
\lim_{r \to \infty} \mu (\phi_{r}A  \, \triangle \, \phi A) 
\ = \ 0  .
\]
Each $\phi \in \Phi$ represents 
an equivalence class $[\phi] = \{ \psi \in \Phi: \mu (\phi=\psi) = 1 \}$ of transformations;
for convenience of notation we identify $[\phi]$ with $\phi$.  
Let $\phi \psi = \phi \circ \psi$ denote the composition of transformations $\phi, \psi \in \Phi$.
For a measure preserving transformation $T$ and $\phi \in \Phi$, let $T_\phi$ denote
the conjugate transformation $\phi^{-1} \circ T \circ \phi$.

\begin{thm} 
\label{counterthm}
Let $T$ be an invertible ergodic measure preserving 
transformation of a Lebesgue probability space $(\X, \calS, \mu)$.
If $\C$ is a class of measurable sets in $\X$ 
such that $H(\C) > 0$, then there exists 
a dense $G_{\delta}$ set $G \subset \Phi$ 
such that for each $\phi \in G$, 
\[
\limsup_{n \to \infty} \,
\sup_{C \in \C} \, \int \Big| \frac{1}{n} \sum_{i=0}^{n-1} I_C(T_\phi^i x) - \mu(C) \Big| \, d\mu 
\ > \ 
0 .
\]
\end{thm}

\vskip.3in

Recall that a measure preserving transformation $T$ of $(\X, \calS, \mu)$ is strong mixing 
(and therefore ergodic) if for all measurable sets $A, B \in \calS$, 
\[
\lim_{n \to \infty} 
\vert \mu (A \cap T^{-n} B) - \mu (A) \mu(B) \vert 
\ = \ 0.
\]
As a corollary of Theorems \ref{mainthm} and \ref{counterthm} we obtain an entropy based
characterization of uniform strong mixing.  An analogous result for weak mixing transformations
is given in Corollary \ref{wkmix} of Section \ref{WkMix}.

\vskip.2in

\begin{cor}
\label{stmix} 
Let $T$ be a strong mixing transformation of $(\X, \calS, \mu)$ 
and let $\C \subseteq \calS$.  

\vskip.1in

\begin{enumerate}
\item[(i)] If $H(\C)=0$ then $T$ is uniform strong mixing on $\C$, in the sense that
\[
\lim_{n \to \infty} \, \sup_{A,B \in \C} 
\left| \mu(A \cap T^{-n}B) - \mu(A) \mu(B) \right| \, = \, 0.
\]

\vskip.1in

\item[(ii)] Suppose that $T$ is invertible.
If $H(\C) > 0$, then there exists a dense $G_{\delta}$ set $G \subset \Phi$ 
such that for each $\phi \in G$,
\[
\limsup_{n \to \infty} \, \sup_{C \in \C} 
\left| \mu (C\cap T_{\phi}^{-n}C) - \mu (C)^2\right| \, > \, 0.
\]
\end{enumerate}
\end{cor}

\vskip.1in

Note that if $T \in \Phi$ is strong mixing then $T_{\phi}$ is strong mixing as well.
As a counterpoint to part (ii) of Corollary \ref{stmix}, we show that for positive entropy
families uniform mixing need not fail for {\em every} strong mixing transformation.
Suppose $S$ and $T$ are invertible commuting transformations on $(\X,\calS,\mu)$ 
that form a strong mixing $\Z^2$-action in the sense that
\[
\mu (T^nS^mA\cap B) \to \mu (A)\mu (B)
\]
for all $A,B\in \calS$ as $|m|+|n|\to \infty$. 
(See \cite{AdSil99,FieFri86,Rud78} for explicit constructions of strong mixing $\Z^2$-actions.)
Let $C\in \calS$ be such that $\mu (C)=\frac12$. 
Define $C_i=S^iC$ and let $\C=\{C_i : i \in \N \}$. 
Since $S$ is strong mixing, $H(\C) = \log{2} > 0$. 
On the other hand,
\[
\mu (T^nC_i\cap C_j) \, = \, \mu (T^nS^iC\cap S^jC) \, = \, \mu (T^nS^{i-j}C\cap C) \, \to \, \mu (C)^2
\]
as $|n|+|i-j|\to \infty$.  Thus 
\[
\lim_{n\to \infty}\sup_{i,j}|\mu (T^nC_i\cap C_j)-\mu (C_i)\mu (C_j)| \ = \ 0
\]
so that $T$ is uniformly mixing over the collection $\C$.

\vskip.2in

\subsection{Overview of the Paper}

The next section is devoted to several elementary approximation properties of entropy.
Section \ref{Pf-thm12} contains the proofs of Theorems \ref{mainthm} and \ref{counterthm}.
The proof of Corollary \ref{stmix} and the statement of an analogous result for weak
mixing are given in Section \ref{StMix}.

\section{Basic Properties of Entropy}
\label{EntSec}

Let $(\X, \calS, \mu)$ be a probability space.  (The Lebesgue assumption is not necessary for
the results of this section.)  
If $\Pi_0$, $\Pi_1$ are finite measurable partitions of $\X$,
the conditional entropy of $\Pi_1$ given $\Pi_0$ is defined by
\[
H_\mu(\Pi_1 \, | \, \Pi_0) 
\, = \,
- \sum_{A \in \Pi_0} \mu(A) \, 
\sum_{B \in \Pi_1}
\frac{\mu (B \cap A)}{\mu (A)} 
\log{\frac{\mu (B \cap A)}{\mu (A)}} .
\]
For $C \in \calS$ define $H_\mu(C \, | \, \Pi) = H_\mu(\{C, C^c\} \, | \, \Pi)$. 
If ${\mathcal D}$ is a collection of measurable sets, 
define $\cup \, {\mathcal D} = \cup_{A \in {\mathcal D}} \, A$.
In the remainder of this section the measure $\mu$ is fixed: in order to reduce notation
we will drop the subscript from the entropy.

The importance of entropy in our analysis derives from its connection 
with the approximation of sets by finite partitions. We present two elementary approximation results, and an alternative
characterization of zero entropy families, whose proofs
are provided for completeness.

\begin{lem} 
\label{entlem1}
For every $\epsilon > 0$, there exists $\delta > 0$ such that 
for every finite measurable partition $\Pi$ and every measurable set $C$ the relation  
$H(C \, | \, \Pi) > \epsilon$ implies
\[ 
\mu 
\left( \, \cup \, \{A \in \Pi: \delta \mu(A) \, \leq \, \mu(A \cap C) \, \leq \, (1-\delta) \mu (A) \} \, \right) 
\, > \, \delta .
\]
\end{lem}

\noindent
{\bf Proof}:
Fix $\epsilon > 0$ and let $0 < \delta < \epsilon / 2$ be such that 
$-x \, \log{x} < \epsilon /4$ when $0 < x < \delta$ or 
$1-\delta < x < 1$. 
Let $\Pi$ be a finite measurable partition and let $C$ 
be a measurable set such that $H(C \, | \, \Pi) > \epsilon$. 
Define the collections
\[
\Pi_1 = \{ A \in \Pi: \mu (A \cap C) < \delta \mu(A) \} 
\ \mbox{ and } \  
\Pi_2 = \{ A \in \Pi: \mu (A \cap C) > (1 - \delta) \mu (A) \} .
\]
Let $\Pi_3 = \Pi \setminus (\Pi_1 \cup \Pi_2)$. 
The definition of $\delta$ and the fact that $H( \{ B, B^c \}) \leq 1$ for
any measure $\mu$ together imply that
\begin{eqnarray*}
H(C \, | \, \Pi) 
& = & 
-\sum_{A \in \Pi_3} \mu (A)
\left[ \frac{\mu (A\cap C)}{\mu (A)} 
\log{\frac{\mu (A\cap C)}{\mu (A)}} + 
\frac{\mu (A\cap C^c)}{\mu (A)} 
\log{\frac{\mu (A\cap C^c)}{\mu (A)}} \right] \\[.1in] 
&  & 
- \sum_{A \in \Pi_3^c} \mu (A)
\left[ \frac{\mu (A\cap C)}{\mu (A)} 
\log{\frac{\mu (A\cap C)}{\mu (A)}} + 
\frac{\mu (A\cap C^c)}{\mu (A)} 
\log{\frac{\mu (A\cap C^c)}{\mu (A)}} \right] \\[.1in] 
& < & 
\mu (\cup \, \Pi_3) + \frac{\epsilon}{2} \mu (\cup \, \Pi_3^c)
\ \leq \ 
\mu (\cup \, \Pi_3) + \frac{\epsilon}{2} .
\end{eqnarray*}
As $H(C \, | \, \Pi) > \epsilon$, it follows that 
$\mu (\cup \, \Pi_3) > \epsilon / 2 > \delta$. 
$\Box$ 

\vskip.3in

\begin{lem}
\label{entlem2}
For every $\epsilon > 0$, there exists $\delta > 0$ such that 
for every finite measurable partition $\Pi$ and every measurable set $C$ satisfying 
$H(C \, | \, \Pi) < \delta$ there is a subcollection $\Pi_0 \subset \Pi$ depending on $C$ 
such that $\mu ( (\cup \, \Pi_0) \, \triangle \, C) < \epsilon$. 
\end{lem}

\noindent
{\bf Proof:} Fix $\epsilon >0$ and let 
$\alpha = - (\epsilon / 4) \log (\epsilon / 4) - 
(1 - \epsilon / 4) \log{(1 - \epsilon / 4)}$. 
Choose $0 < \delta < \epsilon \, \alpha / 4$, and let $\Pi$ and $C$ be such that
$H(C \, | \, \Pi) < \delta$.  Define the collections
\[
\Pi_0 = \{ A \in \Pi: \mu (A \cap C) > (1 - \epsilon) \mu (A) / 4 \} 
\ \mbox{ and } \  
\Pi_1 = \{ A \in \Pi: \mu (A \cap C) < \epsilon \mu(A) / 4 \} .
\]
Let $\Pi_2 = \Pi \setminus (\Pi_0 \cup \Pi_1)$ and define $B_i = \cup \, \Pi_i$.
Note that 
\[
H(C \, | \, \Pi)
\geq 
\sum_{A \in \Pi_2} \mu (A) \, \alpha 
 = 
\alpha \, \mu (B_2) ,
\] 
and therefore $\mu (B_2) < \delta / \alpha < \epsilon / 4$. 
The definition of $\Pi_1$ ensures that $\mu (C \cap B_1) < \epsilon / 4$, so that
\[
\mu (C \setminus B_0)
\, = \,
\mu (C \cap B_1) + \mu (C \cap B_2) 
\, < \,
\frac{\epsilon}{4} + \frac{\epsilon}{4} 
\, = \,
\frac{\epsilon}{2} .
\]  
Finally, the definition of $\Pi_0$ ensures that 
$\mu (B_0 \setminus C) < \epsilon / 4$, and the result follows.
$\Box$

\vskip.3in

\begin{lem}
\label{entlem3}
Let $\C \subseteq \calS$ be a family of measurable subsets of $X$.  
The entropy $H(\C)=0$ if and only if for each $\delta > 0$ there exists 
a finite measurable partition $\Pi$ such that 
$H(C \, | \, \Pi) \leq \delta$ for all $C \in \C$. 
\end{lem}

\noindent
{\bf Proof:} The proof follows from standard facts about entropy.  
Suppose that the approximation condition in the statement of the lemma holds. 
Let $\delta > 0$ and let $\mathbb{S} = (C_1, C_2, \ldots)$ 
be an ordered sequence of sets in $\C$.  For each $n \geq 1$ define 
$\Pi_n = \vee_{i=1}^n C_i$. 
By standard results,
\[ 
H(\Pi_n | \Pi) \, \leq \, \sum_{i=1}^{n} H(C_i \, | \, \Pi)
\ \ \mbox{ and } \ \ 
H(\Pi_n \, | \, \Pi) = H(\Pi_n \vee \Pi) - H(\Pi).
\]
Let $\Pi$ be a finite partition such that 
$H(C \, | \, \Pi) < \delta$ for all $C \in \C$.  Then 
\begin{eqnarray*}
\frac{1}{n} H(\vee_{i=1}^n C_i)
& = & 
\frac{1}{n} H(\Pi_n)
\ \leq \ 
\frac{1}{n} H(\Pi_n \vee \Pi) \ = \ \frac{1}{n} [ H(\Pi_n \, | \, \Pi)+H(\Pi) ] \\ 
& \leq & 
\frac{1}{n} \sum_{i=1}^n H(C_i \, | \, \Pi) + \frac{H(\Pi)}{n}
\ \leq \ \delta + \frac{H(\Pi)}{n}
\ < \ 
2 \delta
\end{eqnarray*}
when $n$ is sufficiently large. Since $\delta$ and $\mathbb{S}$ were arbitrarily,
$H(\C) = 0$. 

Now suppose that the approximation condition fails to hold. Then there exists 
$\delta > 0$ such that for any finite partition $\Pi$ there is a set $C \in \C$ 
such that $H(C \,| \, \Pi) \geq \delta$.   Let $C_1$ be any set in $\C$.  For 
each $i \geq 2$ choose $C_i \in \C$ such that $H(C_i, | \, \Pi_{i-1}) \geq \delta$,
where $\Pi_r =  \vee_{j=1}^r C_j$. 
By the chain rule for entropy, 
\[
H(\Pi_n) \ = \ \sum_{i=2}^n H(C_i \, | \, \Pi_{i-1}) + H(\Pi_1) \ \geq \ (n-1) \delta .
\]
Thus $H(\mathbb{S}) \geq \delta$ where $\mathbb{S} = (C_1, C_2, \ldots)$, and 
it follows that $H(\C) \geq \delta$.  $\Box$

\vskip.2in

\noindent
{\bf Definition:} 
Fix $k \geq 1$ and let $\calS^k$ be $k$-fold Cartesian product of the sigma field $\calS$.
A sequence of functions $f_n: \calS^k \to \real$, $n \geq 1$, is uniformly continuous with respect to $\mu$ if 
for each $\epsilon > 0$ there exists a $\delta > 0$ 
such that $A_1, \ldots, A_k, B_1,\ldots, B_k \in \calS$ and 
$\mu (A_i \triangle B_i) < \delta$ for $i=1, 2, \ldots, k$ imply 
$\vert f_n (A_1, \ldots, A_k) - f_n (B_1, \ldots, B_k) \vert < \epsilon$ for all $n \geq 1$. 

\vskip.2in


\begin{prop}
\label{UCS}
Let the functions $f_n: \calS^k \to \real$, $n \geq 1$, be a uniformly continuous with respect to $\mu$ and such that 
$\lim_{n \to \infty} f_n(A_1,\ldots,A_k) = 0$ for all $A_1, \ldots, A_k \in \calS$. 
If $\C \subseteq \calS$ is a class of measurable sets such that 
$H_\mu(\C)=0$, then 
\[
\lim_{n \to \infty} \, \sup_{C_1, \ldots, C_k \in \C} \, 
f_n(C_1, \ldots, C_k) 
\ = \ 0.
\] 
\end{prop}

\noindent
{\bf Proof}: Suppose that $H(\C) = 0$.  Fix $\epsilon > 0$ and let $\delta > 0$ be as
in the definition of uniform continuity above.  Combining Lemma \ref{entlem3} and 
Lemma \ref{entlem2}, we see that there exists a finite measurable partition $\Pi$ 
such that for every $C \in \C$ there is an approximating set $\hat{C} \in \calS(\Pi)$
such that $\mu(C \, \triangle \, \hat{C}) < \delta$.
Here $\calS(\Pi)$ is the (finite) family of all finite unions of cells of $\Pi$.  
Let $N \geq 1$ be so large that $| f_n( B_1, \ldots, B_k) | \leq \epsilon$ for
every $n \geq N$ and every $B_1, \ldots, B_k \in \calS(\Pi)$.  Then for
each $n \geq N$, and each sequence of sets $C_1, C_2, \ldots, C_k \in \C$, 
\[
\vert f_n (C_1, \ldots, C_k) \vert  
\ \leq \ 
\vert f_n (C_1, \ldots, C_k) - f_n (\hat{C}_1, \ldots, \hat{C}_k) \vert  
\ + \  
\vert f_n (\hat{C}_1, \ldots, \hat{C}_k) \vert 
\ \leq \  
2 \epsilon .
\]
As $\epsilon > 0$ and the sets $C_1,\ldots, C_k$ were arbitrary, the proof is complete.  $\Box$

\vskip.3in

\section{Proof of Theorems \ref{mainthm} and \ref{counterthm}}
\label{Pf-thm12}

\subsection{Sufficiency of Zero Entropy}

\begin{prop}
\label{zentume}
Let $(\X, \calS, \mu)$ be a probability space, and 
let $\C \subseteq \calS$. 
If $H_\mu(\C) = 0$ then 
for every ergodic measure preserving transformation $T$ of $(X,\calS,\mu)$,
\[
\sup_{C \in \C} \int \left| \, \frac{1}{n} \sum_{i = 0}^{n-1} I_C(T^i x) - \mu(C) \, \right| \, d\mu
\ \to \ 0  
\]
as $n$ tends to infinity.
\end{prop}

\noindent
{\bf Proof:} For each $n \geq 1$ define the set function $f_n: \calS \to \real$ by 
\[
f_n(A) \ = \ \int \left| \, \frac{1}{n} \sum_{i = 0}^{n-1} I_A(T^i x) - \mu(A) \, \right|  d\mu . 
\]
The mean ergodic theorem ensures that $f_n(A) \to 0$ for each $A \in \calS$.
Let $A, B$ be measurable sets. By standard arguments,
\begin{eqnarray*}
f_n(A) - f_n(B) 
& = & 
\int
\left( \Big| \frac{1}{n} \sum_{i = 0}^{n-1} I_A(T^i x) -  \mu (A) \Big|
- \Big|\frac{1}{n} \sum_{i = 0}^{n-1} I_B(T^i x)  - \mu(B)  \Big| \right) d\mu \\[.1in]
& \leq & 
\int
\left|  \left( \frac{1}{n} \sum_{i = 0}^{n-1} I_A(T^i x) -  \frac{1}{n} \sum_{i = 0}^{n-1} I_B(T^i x) \right)
+ \mu(B) - \mu(A) \, \right| d\mu \\[.1in]
& \leq &
\frac{1}{n} \sum_{i = 0}^{n-1} \int_X | I_A(T^i x) - I_B(T^i x) | \, d\mu
\, + \,
| \mu(B) - \mu(A) | \\[.1in]
& \leq &
2 \mu(A \, \Delta \, B) .
\end{eqnarray*}
Thus $\{ f_n \}$ is uniformly continuous with respect to $\mu$, and  
the conclusion of the proposition follows from Proposition \ref{UCS}. $\Box$

\vskip.2in

\subsection{Strong Necessity of Zero Entropy}


Consider for the moment the special case in which $\X = [0,1)$ is the unit interval
equipped with its Borel subsets ${\mathcal B}$, and $\mu$ is Lebesgue measure.  
As above, let $\Phi$ be the set of invertible measure preserving transformations 
of $([0,1), {\mathcal B}, \mu)$ endowed with the weak topology.  
Let $E_1, E_2, \ldots $ be a countable sequence of measurable sets generating
${\mathcal B}$.  Define the distance between transformations $\phi, \psi \in \Phi$ by 
\[
d(\phi,\psi) \ = \ 
\sum_{i=1}^{\infty}\frac{1}{2^i}
\left[ \, \mu ( \phi E_i \, \triangle \, \psi E_i ) 
\, + \,  \mu ( \phi^{-1} E_i \, \triangle \, \psi^{-1} E_i ) \, \right].
\]
It follows from standard results \cite{Hal56} that $(\Phi,d)$ is a complete metric space,
and that the topology generated by $d(.,.)$ coincides with the weak topology on $\Phi$. 

\vskip.1in

\begin{prop}
\label{denseprop} 
Let $T$ be an invertible ergodic measure preserving transformation of the space
$([0,1), {\mathcal B}, \mu )$
and let $\C \subseteq {\mathcal B}$ be a family of Borel sets.
If $H(\C) > 0$, then there exists $\delta > 0$ such that 
\[
G_N (T, \delta) 
\ = \ 
\left\{ \phi \in \Phi : \sup_{n \geq N} \, \sup_{C \in \C} 
\int \Big\vert \frac{1}{n} \sum_{i=0}^{n-1} I_C(T_{\phi} (x)) - \mu (C) \Big\vert \, d\mu(x)  > \delta 
\right\}.
\] 
is dense in $\Phi$ for every $N \geq 1$.
\end{prop}

\noindent
{\bf Proof:}  By Lemma \ref{entlem3}, there exists $\delta_0 > 0$ such that, 
for every finite measurable partition $\Pi$ of $[0,1)$ there is some set $C \in \C$ 
such that $H(C \,| \, \pi) > \delta_0$.   It then follows from Lemma \ref{entlem1} 
that for every finite measurable partition $\Pi$
\[
\mu \big( \cup \{ A \in \Pi : \delta_{1} \, \mu(A) \, \leq \, \mu(A \cap C) \, \leq \, (1-\delta_{1}) \, \mu (A) \} \, \big) 
\ > \ \delta_{1} 
\]
for some constant $\delta_{1} > 0$ depending only on $\delta_0$.  
Define $\delta =\delta_1^3 / 64$.
Fix $N \geq 1$ and let $\epsilon > 0$ and $\phi \in \Phi$ be given.  
We will construct invertible measure preserving transformations 
$S$, $\psi$ and $\tau$ such that the following three conditions hold
\[ 
S_{\psi} = T_{\phi \tau^{-1} \psi}, \ \ 
d(\phi, \phi \tau^{-1} \psi ) < \epsilon, \ \
\mbox{and} \ 
\sup_{n \geq N} \, \sup_{C \in \C} 
\int \left\vert \frac{1}{2n}\sum_{i=0}^{2n-1} 
I_C(S_{\psi}^i x) - \mu (C) \right\vert d\mu 
\ > \ \delta .
\]
To this end, choose integers $n \geq N$ and $m \geq 1$ such that 
\be
\label{choice}
\frac{16}{n} \, < \, \epsilon 
\ \ \mbox{ and } \ \   
\sum_{i=m}^{\infty}\frac{1}{2^i}=\frac{1}{2^{m-1}}<\frac{\epsilon}{2}.
\ee

Let $D_0, D_1,\ldots ,D_{2n-1} \subseteq [0,1)$ be a Rohklin tower of height $2n$ for $T_\phi$: in other words,
the sets $D_i$ are disjoint and satisfy the relations
\be
\label{Rohklin}
D_i =T_{\phi}^iD_0
\ \ 
0 \leq i \leq 2n-1
\ \mbox{ and } \  
\mu(\bigcup_{i=0}^{2n-1}T_{\phi}^iD_0) > 1-\frac{1}{2n} .
\ee
Let $E_1, E_2,\ldots$ be the sequence of generating sets used to define the metric $d(.,.)$ on $\Phi$. 
With $m \geq 1$ as in (\ref{choice}) define finite partitions
$
\Gamma_0 = [ \bigvee_{i=1}^{m} E_i ] \vee [ \bigvee_{i=1}^{m} \phi^{-1}E_i ]
$ 
and
$
\Gamma_1 = \bigvee_{i=0}^{2n-1} T_{\phi}^{-i} \Gamma_0 
$
of $[0,1)$.  Let
\[
\Pi_0 \ = \ \{ A \cap D_0 : A \in \Gamma_1 \} \ = \ \{ A_0,\ldots, A_K \}
\]
be the partition of $D_0$ induced by $\Gamma_1$, and let 
\[
\Pi_1 \ = \  \{ T_{\phi}^i (A_k): A_k \in \Pi_0, 0 \leq i < 2n \}
\]
be a partition of the levels of the Rohklin tower obtained from the partition of the base
$D_0$ using the map $T_{\phi}$.  Note that $\tilde{\Pi}_1 = \Pi_1 \cup \{ (\cup_i D_i)^c \}$
is a partition of $[0,1)$.  Our choice of $\Gamma_1$ ensures that
for each $A \in \Pi_1$ and each $1 \leq j \leq m$ either $A \subset E_j$ or $A \subset E_j^c$, 
and either $A \subset \phi^{-1}E_j$ or $A \subset \phi^{-1}E_j^c$. 
Define sets
\be
\label{fjdef}
F_j = \cup \{ A \in \Pi_1 : A \subset E_j \} \ \mbox{ and } \ 
F_j^{\prime} = \cup \{ A \in \Pi_1 : A \subset \phi^{-1} (E_j) \} .
\ee 
It follows from (\ref{Rohklin}) that $\mu (F_j \, \triangle \, E_j) < \frac{1}{2n}$ 
and $\mu (F_j^{\prime} \, \triangle \, \phi^{-1}E_j) < \frac{1}{2n}$ 
for each $1\leq j \leq m$.

Choose $C \in \C$ such that $H(C \, | \, \tilde{\Pi}_1) > \delta$. 
Let $\beta_0 = 0$ and define $\beta_k =  \sum_{l=0}^{k-1} \mu (A_l)$
for $1 \leq k \leq K$.  
Note that for each $0 \leq k \leq K$ the set $A_k$ is isomorphic to
the interval $[\beta_k, \beta_{k+1})$.  For each $k$ define an invertible measure 
preserving map
$
\psi_k : A_k \to [\beta_k, \beta_{k+1})
$
such that 
\[
\psi_k (C \cap A_k) =
[\beta_k, \beta_k + \mu (C \cap A_k) )
\]
Let $\alpha_0 = 0$ and define $\alpha_i = \sum_{l=0}^{i-1} \mu (D_i) = i \, \mu(D_0)$
for $1 \leq i \leq 2n-1$.  Define $\psi: D_0 \to [0, \alpha_1)$ such that 
$\psi (x) = \psi_k (x)$ for $x \in A_k$. 
Extend $\psi$ to an invertible measure preserving map on $[0,1)$ as follows: 
\begin{eqnarray*}
\psi : D_i & \to & [\alpha_i, \alpha_{i+1}) \\[.1in] 
\psi : T_{\phi}^i (A_k) & \to & [\alpha_i + \beta_k, \alpha_i + \beta_{k+1}) ) \\[.1in] 
\psi : C\cap T_{\phi}^i (A_k) & \to & [\alpha_i + \beta_k, \alpha_i + \beta_{k} + \mu (C \cap T_{\phi}^i A_k) ) \\[.1in]
\psi : (\cup_{i} D_i)^c & \to & [\alpha_{2n-1},1).
\end{eqnarray*}
Thus $\psi$ maps each set $T_\phi^i(A_k) \in \Pi_1$ to the interval 
$[\alpha_i + \beta_k, \alpha_i + \beta_{k+1})$, and maps the 
intersection of $C$ and $T_\phi^i(A_k)$ to the left side of the interval.
Define $S: [0,1) \to [0,1)$ such that 
\begin{eqnarray*} 
S(x)= 
\left\{\begin{array}{ll}
x + \mu (D_0) & \ \mbox{if $\ x\in [0, \alpha_{2n-2})$} \\ 
\psi \circ T_{\phi} \circ \psi^{-1} (x) & \ \mbox{if $\ T_{\phi}(\psi^{-1} x ) \in (\bigcup_{i=0}^{2n-1} D_i)^c $} \\ 
S^{-2n+1} \circ \psi \circ T_{\phi}^{2n} \circ \psi^{-1} (x) & \ \mbox{if \ $T_{\phi}(\psi^{-1} x )\in D_0$}
\end{array}
\right.
\end{eqnarray*}
One may readily verify that $S$ is invertible and measure preserving.  
Define $\tau: [0,1) \to [0,1)$ such that 
\begin{eqnarray*} 
\tau (x) = 
\left\{\begin{array}{ll}
\psi (x) & \ \mbox{if } x \in (\bigcup_{i=0}^{2n-2}D_i)^c \\ 
(S^{-i} \circ \psi \circ T_{\phi}^i)(x) & \ \mbox{if } x \in D_{2n-i-1}, \,1 \leq i < 2n
\end{array}
\right.
\end{eqnarray*}
It is easy to show by induction that $\tau$ is well defined on $[0,1)$. 


The transformation $\psi$ is defined to map small pieces of the set $C$ to the left side 
of the subcolumns determined by the $\beta_k$'s.  This constructs a new tower, 
consisting primarily of intervals, such that the transformation $S$ has a poor average 
on $\psi (C)$, as it maps points up the tower.  The map $S$ is defined on the top 
portion of the tower to be isomorphic to $T_{\phi}$. 
Thus, the map $S_{\psi}$ has a poor average on $C$, and is isomorphic to $T_{\phi}$. 
Also, $\psi$ has been defined such that $S_{\psi}$ approximates $T_{\phi}$ 
in the weak topology. 
The map $\tau$ is the isomorphism used to map $S$ to $T_{\phi}$. 
In particular, $S_{\psi} = T_{\phi \tau^{-1} \psi}$, and 
$\phi \tau^{-1} \psi$ is close to $\phi$ in our metric. 
We give a precise proof of these facts below.

\vskip.2in

\noindent
{\bf Claim 1:} $S_{\tau} = T_{\phi}$. 

\vskip.1in

\noindent
{\bf Proof of Claim 1:} 
We show inductively that 
$
\tau^{-1} \circ S \circ \tau = T_{\phi} \mbox{ on } D_{2n-i-1} \mbox{ for } 1 \leq i < 2n. 
$
Let $i=1$.  Then for $x \in D_{2n-2}$,
\[
(\tau^{-1} \circ S \circ \tau) (x) = 
(\tau^{-1}\circ S \circ S^{-1} \circ \tau) (T_{\phi} \, x)
= T_{\phi} \, x.
\] 
Suppose now that $\tau^{-1}\circ S \circ \tau = T_{\phi}$ for 
$x \in D_{2n-i-1}$, and consider $x \in D_{2n-i-2}$.  The definition of $\tau$ 
and the inductive hypothesis ensure that
\begin{eqnarray*}
(\tau^{-1} \circ S \circ \tau) (x)
& = &
(\tau^{-1} \circ S \circ S^{-i-1} \circ \tau) (T_{\phi}^{i+1}x) \\
& = &
(\tau^{-1} \circ S^{-i} \circ \tau) (T_{\phi}^{i+1}x) \\
& = & 
(T_{\phi}^{-i} \circ T_{\phi}^{i+1}) (x) \\
& = & 
T_{\phi} \, x.
\end{eqnarray*}
Suppose $x \in (\bigcup_{i=0}^{2n-2}D_i)^c$. Either 
$T_{\phi}(x ) \in (\bigcup_{i=0}^{2n-1} D_i)^c$ or 
$T_{\phi}(x ) \in D_0$. 
If $T_{\phi}(x ) \in (\bigcup_{i=0}^{2n-1} D_i)^c$, then 
\begin{eqnarray*}
(\tau^{-1} \circ S \circ \tau) (x)
& = &
(\tau^{-1} \circ S ) (\psi x) \\ 
& = & \tau^{-1} \circ \psi \circ T_{\phi} \circ \psi^{-1} (\psi x) \\ 
& = & \tau^{-1} (\psi (T_{\phi}x)) = T_{\phi}x . 
\end{eqnarray*}
If $T_{\phi}(x) \in D_0$, then 
\begin{eqnarray*}
(\tau^{-1} \circ S \circ \tau) (x)
& = &
(\tau^{-1} \circ S ) (\psi x) \\ 
& = & \tau^{-1} \circ S^{-2n+1} \circ \psi \circ T_{\phi}^{2n} \circ \psi^{-1} (\psi x) \\ 
& = & \tau^{-1} \circ S^{-2n+1} \circ \psi \circ T_{\phi}^{2n} (x) \\ 
& = & (T_{\phi}^{-2n+1} \circ \psi^{-1} \circ S^{2n-1}) (S^{-2n+1} \circ \psi \circ T_{\phi}^{2n})(x) \\ 
& = & T_{\phi} (x) . \Box 
\end{eqnarray*}

\vskip.2in 

The transformation isomorphic to $T_{\phi}$ with the required "bad" average is $S_{\psi}$. 
The previous claim implies $S_{\psi} = (T_{\phi})_{\tau^{-1} \psi} = T_{\phi \tau^{-1} \psi}$. 
Hence, we wish to show that $\phi \circ \tau^{-1} \circ \psi$ is a small perturbation of $\phi$ 
in our metric space. 

\vskip .2in 

\noindent
{\bf Claim 2:} $d(\phi,\phi \circ \tau^{-1} \circ \psi ) < \epsilon$. 

\vskip.1in

\noindent
{\bf Proof of Claim 2:}  For $x \in \bigcup_{i=0}^{2n-1} D_i$ and $0 \leq k \leq K$
the definitions of $\tau$ and $S$ lead to the following chain of equivalences:
\begin{eqnarray*}
x \in T_{\phi}^i (A_k)
& \Leftrightarrow & 
\psi (x) \in [\alpha_i +\beta_k, \alpha_i +\beta_{k+1}) \\
& \Leftrightarrow & 
S^{2n-i-1}(\psi (x)) \in [\alpha_{2n-1} + \beta_k, \alpha_{2n-1} +\beta_{k+1}) \\
& \Leftrightarrow & 
\psi^{-1} (S^{2n-i-1} (\psi (x))) \in T_{\phi}^{2n-1} (A_k)  \\ 
& \Leftrightarrow &
T_{\phi}^{-2n+i+1}(\psi^{-1} (S^{2n-i-1} (\psi (x)))) \in T_{\phi}^{i} (A_k)
\end{eqnarray*} 
Therefore
\be
\label{spsi1}
\tau^{-1} \circ \psi (T_{\phi}^i A_k) = T_{\phi}^i A_k \ \mbox{ for } \ 0 \leq i \leq 2n-1 .
\ee
This implies that 
\be
\label{spsi2}
\psi^{-1}\circ \tau (T_{\phi}^i A_k) = T_{\phi}^i A_k \ \mbox{ for } \ 0 \leq i \leq 2n-1 .
\ee
Let $F_j = \cup \{ A \in \Pi_1 : A \subset E_j \}$ be the $\Pi_1$-approximation of $E_j$ 
defined above.  By the triangle inequality, for $1 \leq j \leq m$,
\begin{eqnarray}
\lefteqn{\mu (\phi (E_j) \, \triangle \, \phi \circ \tau^{-1}\circ \psi (E_j))} \nonumber \\[.1in] 
& \leq & 
2\mu (E_j \, \triangle \, F_j) \, + \, 
\mu (F_j \, \triangle \, \tau^{-1} \circ \psi (F_j)) \label{triangle}. 
\end{eqnarray}
It was shown above that $\mu (E_j \, \triangle \, F_j) < (2n)^{-1}$.  Moreover, (\ref{spsi1}) implies that 
$$\mu (F_j \, \triangle \, \tau^{-1} \circ \psi (F_j)) = 0$$ 
and hence 
$\mu (\phi (E_j) \, \triangle \, \phi \circ \tau^{-1}\circ \psi (E_j)) < 4/n$ for $1 \leq j \leq m$.
A similar argument using $F_j^{\prime}$ in place of $F_j$, and applying (\ref{spsi2}) shows that
$\mu (\phi^{-1}(E_j) \, \triangle \, \psi^{-1}\circ \tau \circ \phi^{-1} (E_j)) < 4/n$ for $1 \leq j \leq m$, and therefore
\[
\sum_{j=1}^{m} 2^{-j} 
[\, \mu (\phi (E_j) \, \triangle \, \phi \circ \tau^{-1} \circ \psi (E_j)) + 
\mu (\phi^{-1} (E_j) \, \triangle \, \psi^{-1} \circ \tau \circ \phi^{-1}(E_j)) \, ] \ < \ \frac{8}{n}.
\]
It follows from the choice of $m$ that $d(\phi, \phi \circ \tau^{-1} \circ \psi) < \epsilon$, as desired.

\vskip.2in

\noindent
{\bf Claim 3:} $\phi \circ \tau^{-1} \circ \psi \in G_n(T,\delta)$. 

\vskip.1in

\noindent
{\bf Proof of Claim 3:}  
Note that $T_{\phi \tau^{-1} \psi} = S_{\psi}$. 
Define $\Delta = \cup \{ A \in \tilde{\Pi}_1: \delta_{1} \mu(A) \leq \mu(A \cap C) \leq (1-\delta_{1}) \mu (A)\}$,
and let $B = \bigcup_{i=0}^{2n-1} S_{\psi}^{-i} D_{2n-1}$.  Our choice of $C$ ensures that
\begin{eqnarray}
\delta_1 \ < \ \int_{[0,1)} I_{\Delta}(x) \, d\mu &=& 
\int_{B} I_{\Delta}(x) \, d\mu 
+ \int_{B^c}I_{\Delta}(x) \, d\mu \\[.05in]
& \leq & \sum_{i=0}^{2n-1} \int_{S_{\psi}^{-i}D_{2n-1}}I_{\Delta}(x)d\mu 
+ \frac{1}{2n} \\[.1in] 
& = & 
\int_{D_0} \sum_{i=0}^{2n-1} I_{\Delta}(S_{\psi}^i x) \, d\mu 
+ \frac{1}{2n}. 
\end{eqnarray}
The inequality above follows from the fact that $\mu(B^c) \leq (2n)^{-1}$; the final expression
follows from a standard change of variables.
Define sets
\[
\Delta_1 \ = \ \bigcup_{k=1}^{K} \psi^{-1} ([\beta_k, \beta_k + \delta_1 \mu (A_k))) \ \subset \ D_0
\]
and 
\[
\Delta_2 \ = \ \bigcup_{k=0}^{K} \psi^{-1} ([\beta_k + (1-\delta_1) \mu (A_k), \beta_{k+1}))  \ \subset \ D_0.
\]
Thus a point $x$ is contained in $\Delta_1$ if (and only if) it is contained in a cell $A_k \in \Pi_0$
and is mapped by $\psi$ to the left $\delta_1$ fraction of the interval 
$[\beta_k, \beta_{k+1})$ associated with $A_k$.  The set $\Delta_2$ has a similar interpretation.
Note that $\mu (\Delta_1)=\mu (\Delta_2)=\delta_1\mu (D_0) \leq \delta_1/ 2n$. 

Recall that $\psi$ maps the intersection of $C$ and the set $T_{\phi}^i A_k$ into the left part of 
the interval $[\alpha_i + \beta_k, \alpha_i + \beta_{k+1})$.  Thus if $S^i(\psi (x))\in \psi (C)$
for some $x \in A_k \cap \Delta_2$, then $S^i (\psi (y)) \in \psi (C)$
for every $y \in A_k \cap \Delta_1$.  It follows that
\[ 
\int_{A_k \cap \Delta_1} I_C (S_{\psi}^i x) \, d\mu 
\ \geq \
\int_{A_k \cap \Delta_2} I_C (S_{\psi}^i x) \, d\mu .
\] 
for $1 \leq k \leq K$, and therefore
\be
\int_{\Delta_1} I_C (S_{\psi}^i x) \, I_{\Delta^c} (S_{\psi}^ix) \, d\mu 
\ \geq \ 
\int_{\Delta_2} I_C (S_{\psi}^i x) \, I_{\Delta^c} (S_{\psi}^i x) \, d\mu.
\ee
Note that if $x \in \Delta_2$ and $S_{\psi}^i(x) \in \Delta$, then $S_{\psi}^i(x) \notin C$.  Thus
\be
\int_{\Delta_2} I_C(S_{\psi}^ix) \, I_{\Delta}(S_{\psi}^ix) \, d\mu \ = \ 0.
\ee
Similarly, if $x \in \Delta_1$ and $S_{\psi}^i(x) \in \Delta_1$, then $S_{\psi}^i(x) \in C$,
and therefore
\be
\int_{\Delta_2} I_C(S_{\psi}^ix) \, I_{\Delta}(S_{\psi}^ix) \, d\mu \ = \ 0.
\ee
Together, the last three displays imply that for each $0 \leq i \leq 2n - 1$, 
\begin{eqnarray}
\lefteqn{ \int_{\Delta_1} I_C (S_{\psi}^i x) \, d\mu 
\ - \ 
\int_{\Delta_2} I_C (S_{\psi}^i x) \, d\mu}  \nonumber \\[.1in]
& \geq & 
\int_{\Delta_1} I_C (S_{\psi}^i x) \, I_{\Delta} (S_{\psi}^i x) \, d\mu 
\ - \ 
\int_{\Delta_2} I_C (S_{\psi}^i x) \, I_{\Delta} (S_{\psi}^i x) \, d\mu  \nonumber \\[.1in]
& = & 
\int_{\Delta_1}  I_C(S_{\psi}^ix) I_{\Delta}(S_{\psi}^ix) \, d\mu  \nonumber \\[.1in]
&=& 
\label{int1}
\int_{\Delta_1} I_{\Delta} (S_{\psi}^i x) \, d\mu 
\end{eqnarray}
Finally, note that if $x \in \Delta_1 \cap A_k$ and $S_\psi^i x \in \Delta$
then $S_\psi y \in \Delta$ for every $y \in A_k$.  Therefore,
\be
\label{int2}
\int_{\Delta_1} I_{\Delta} (S_{\psi}^i x) \, d\mu 
\ \geq \
\frac{\mu (\Delta_1)}{\mu (D_0)} 
\int_{D_0} I_{\Delta} (S_{\psi}^i x) \, d\mu \\ 
\ > \ 
\delta_1 (\delta_1 - \frac{1}{2n}) / 2n. 
\ee
Averaging over $i = 0, \ldots, 2n - 1$, the inequalities (\ref{int1}) and (\ref{int2})  
imply that either 
\be
\label{bavg1}
\int_{\Delta_1}
\Big| \frac{1}{2n}\sum_{i=0}^{2n-1} I_C(S_{\psi}^ix) - \mu (C) \Big| d\mu 
\ > \ 
\delta_1\frac{(\delta_1 - \frac{1}{2n})}{4n} 
\ee
or 
\be
\label{bavg2}
\int_{\Delta_2} \Big| \frac{1}{2n}\sum_{i=0}^{2n-1} I_C(S_{\psi}^ix) - \mu (C) \big|  d\mu 
\ > \
\delta_1\frac{(\delta_1 - \frac{1}{2n})}{4n}. 
\ee

Assume that inequality (\ref{bavg1}) holds.  One may reason from (\ref{bavg2}) 
in a similar fashion.  By the triangle inequality and an elementary change of variables,
\begin{eqnarray*}
\lefteqn{ \int_{\Delta_1} \Big| \frac{1}{2n}\sum_{i=0}^{2n-1} 
I_C(S_{\psi}^ix)-\mu (C) \Big| d\mu } \\[.05in] 
& \leq & 
\int_{\Delta_1} \Big| \frac{1}{2n} \sum_{i=-1}^{2n-2} I_C(S_{\psi}^i x) - \mu (C) \Big| d\mu 
\, + \, 
\frac{1}{2n} \int_{\Delta_1} \Big| I_C (S_{\psi}^{2n-1} x) - I_C (S_{\psi}^{-1} x) \Big| d\mu \\[.05in]  
& = & \int_{S_{\psi}(\Delta_1)}\Big| \frac{1}{2n}\sum_{i=0}^{2n-1} 
I_C(S_{\psi}^ix)-\mu (C)\Big| d\mu + 
\frac{1}{2n} \int_{\Delta_1} \Big| I_C (S_{\psi}^{2n-1} x) - I_C (S_{\psi}^{-1} x) \Big| d\mu \\[.05in]  
& < & 
\int_{S_{\psi}(\Delta_1)}\Big| \frac{1}{2n}\sum_{i=0}^{2n-1} I_C (S_{\psi}^i x) - \mu (C) \Big| d\mu 
\, + \, 
\frac{\delta_1}{2 n^2}
\end{eqnarray*}
where in the last step we have used the fact that 
$\mu(\Delta_1) \leq \delta_1 \mu(D_0) \leq \delta_1 / n$.
Combining the last display with (\ref{bavg1}) yields the inequality
\[
\int_{S_{\psi}(\Delta_1)} 
\Big| \frac{1}{2n} \sum_{i=0}^{2n-1} I_C(S_{\psi}^i x) - \mu (C) \Big| d\mu 
\ > \
\frac{\delta_1}{4n}(\delta_1 - \frac{1}{2n}) - \frac{\delta_1}{2n^2} \\ 
\ > \ \frac{\delta_1}{4n}(\delta_1 - \frac{3}{n}).
\]
By a similar argument, for $1 \leq i \leq 2n-1$, we have
\[
\int_{S_{\psi}^j(\Delta_1)} 
\Big| \frac{1}{2n}\sum_{i=0}^{2n-1} I_C(S_{\psi}^ix) - \mu (C)\Big| d\mu 
\ > \
\frac{\delta_1}{4n} (\delta_1 - \frac{2j+1}{n}).
\]
Now choose $r$ such that 
$\delta_1 / 8 < r / n < \delta_1 / 4$. 
Then 
\begin{eqnarray*}
\lefteqn{\int_{X}\Big| \frac{1}{2n}\sum_{i=0}^{2n-1} I_C(S_{\psi}^i x) - \mu(C) \Big| d\mu } \\ 
&\geq& 
\int_{\bigcup_{j=0}^{r-1}S_{\psi}^j(\Delta_1)}\Big| \frac{1}{2n}\sum_{i=0}^{2n-1} 
I_C(S_{\psi}^ix)-\mu (C)\Big| d\mu \\ 
&=& 
\sum_{j=0}^{r-1}\int_{S_{\psi}^j(\Delta_1)}\Big| \frac{1}{2n}\sum_{i=0}^{2n-1} 
I_C(S_{\psi}^ix)-\mu (C)\Big| d\mu \\ 
&>&
\sum_{j=0}^{r-1}\frac{\delta_1}{4n}(\delta_1 - \frac{2j+1}{n}) \\ 
&>& r(\frac{\delta_1^2}{8n}) > \frac{\delta_1^3}{64}.\Box
\end{eqnarray*}
It follows that $S_{\psi}\in G_n(T_{\phi},\delta)$.  This establishes Claim 3, and completes the 
proof of the proposition.

\vskip.2in 

\begin{prop}
\label{openprop}
For each positive integer $N$ and $\delta >0$, $G_N(T,\delta)$ 
is open in $\Phi$. 
\end{prop}
\noindent 
{\bf Proof:} Let $\phi \in G_N(T,\delta)$. Then there exists $n \geq N$ and $C \in \C$ 
such that 
\[
\delta' \, := \, \Big| n^{-1} \sum_{i=0}^{n-1} I_C(T_\phi x) -\mu (C) \Big|  > \delta.
\] 
Let $\epsilon > 0$ be so small that for each $\psi \in \Phi$ with $d(\psi , \phi) < \epsilon$ one has 
\[
\mu (\psi (C) \, \triangle \, \phi (C)) < \frac{(\delta^{\prime} - \delta)}{2} 
\ \mbox{ and } \  
\mu ( \psi^{-1} \circ T^{-i} \circ \phi (C) \, \triangle \, \phi^{-1} \circ T^{-i} \circ \phi (C) ) 
< \frac{(\delta^{\prime} - \delta)}{2} 
\]
for $0 \leq i \leq n-1$.  For each such $\psi$ we claim that 
$\mu(T_{\psi}^{-i}(C) \, \triangle \, T_{\phi}^{-i}(C))< \delta^{\prime} - \delta$
for $1 \leq i \leq n-1$.
To see this, note that by the triangle inequality, 
\begin{eqnarray*}
\mu (T_{\psi}^{-i}(C) \, \triangle \, T_{\phi}^{-i}(C)) 
& = & 
\mu (\psi^{-1} \circ T^{-i} \circ \psi (C) \ \triangle \ \phi^{-1} \circ T^{-i} \circ \phi (C)) \\[.1in] 
& \leq & 
\mu (\psi^{-1} \circ T^{-i} \circ \psi (C) \ \triangle \ \psi^{-1} \circ T^{-i} \circ \phi (C)) \\[.1in] 
& + & 
\mu (\psi^{-1} \circ T^{-i} \circ \phi (C) \ \triangle \ \phi^{-1} \circ T^{-i} \circ \phi (C)) \\[.1in] 
& \leq & 
\mu (\psi (C) \ \triangle \ \phi (C)) 
\ + \ 
\frac{(\delta^{\prime} - \delta)}{2} 
\ < \ 
\delta^{\prime} - \delta .
\end{eqnarray*}
\newline 
\noindent 
As a consequence, we have
\begin{eqnarray*}
\lefteqn{ \int \Big| \frac{1}{n}\sum_{i=0}^{n-1} I_C(T_{\psi}^ix) - \mu (C) \Big| \, d\mu} \\[.1in] 
& \geq & 
\int \Big| \frac{1}{n} \sum_{i=0}^{n-1} I_C(T_{\phi}^ix) - \mu (C) \Big| \, d\mu 
- \int \Big| \frac{1}{n}\sum_{i=0}^{n-1} I_C(T_{\psi}^ix) - \frac{1}{n}\sum_{i=0}^{n-1} I_C(T_{\phi}^ix) \Big| \, d\mu \\[.1in] 
& \geq & 
\delta^{\prime} - 
\frac{1}{n} \sum_{i=0}^{n-1} \mu (T_{\psi}^{-i}(C) \, \triangle \, T_{\phi}^{-i}(C)) \\[.1in]
& > & 
\delta^{\prime} - \frac{1}{n} n (\delta^{\prime} - \delta)
\ = \ 
\delta.
\end{eqnarray*}
It follows that $G_N(T,\delta)$ contains the ball $\{ \psi : d(\phi , \psi ) <\epsilon \}$.  This completes the proof.
$\Box$

\vskip.1in

\subsection{Proof of Theorem \ref{counterthm}}

By applying a measure space isomorphism, it is enough to prove the theorem for the case 
where $\X = [0,1)$, $\calS = {\mathcal B}$, $\C \subseteq {\mathcal B}$, and $\mu$ is 
Lebesgue measure.  By Proposition \ref{denseprop} there exists $\delta > 0$, such that 
$G_N (T,\delta)$ is dense for all $N \geq 1$; by Proposition 2 each set $G_N(T,\delta)$ is open in $\Phi$. 
Define $G = \bigcap_{N=1}^{\infty} G_N(T,\delta)$.  As $(\Phi, d(.,.))$ is complete,
it follows from the Baire category theorem that $G$ is dense in $\Phi$.
If $\phi \in G$, then the definition of $G_N(T,\delta)$ ensures that for each $N \geq 1$, there exists
$n \geq N$ and $C_n \in \C$ such that
\[
\int \Big| \frac{1}{n}\sum_{i=0}^{n-1} I_{C_n}(T_{\phi}^i x) - \mu (C_n) \Big| \, d\mu \ > \ \delta
\]
It follows that the uniform mean ergodic theorem for $\C$ fails to hold for each transformation $T_\phi$
with $\phi \in G$.

\subsection{Proof of Theorem \ref{mainthm}}  

Theorem \ref{mainthm} follows immediately from Proposition \ref{zentume} and Theorem \ref{counterthm}.

\vskip.3in

\section{Uniform Strong and Weak Mixing for Zero-Entropy Families}
\label{StMix}

\noindent
{\bf Proof of Corollary \ref{stmix}:} For each $n \geq 1$ define the set function $f_n: \calS^2 \to \real$ by 
$f_n(A,B) = \mu (A \cap T^{-n} B) - \mu (A) \mu (B)$. 
Let $A_1, A_2, B_1, B_2 \in \calS$ be measurable sets. Then by standard arguments,
\begin{eqnarray*} 
f_n(A_1, B_1) - f_n(A_2, B_2)
& \leq & 
\vert \mu (A_1\cap T^{-n}B_1) - \mu (A_2\cap T^{-n}B_2) \vert \\
& & + \
\vert \mu (A_1)\mu (B_1) - \mu (A_2)\mu (B_2) \vert \\ [.04in]
& \leq & 
\mu (A_1\triangle A_2) + \mu (T^{-n}B_1\triangle T^{-n} B_2) \\ 
& & + \  
\vert \mu (A_1)\mu (B_1) - \mu (A_1)\mu (B_2) \vert 
\ + \ 
\vert \mu (A_1)\mu (B_2) - \mu (A_2)\mu (B_2)\vert \\ [.04in]
& \leq & 
\mu (A_1\triangle A_2) \, + \, \mu (T^{-n}(B_1\triangle B_2)) \\ 
& & + \ 
\mu (A_1) \mu (B_1\triangle B_2) \, + \, \mu(B_2) \mu (A_1\triangle A_2) \\ [.04in]
& \leq & 
2 \mu (A_1\triangle A_2) \, + \, 2 \mu (B_1\triangle B_2) .
\end{eqnarray*}
Thus the sequence $\{ f_n \}$ is uniformly continuous.  Part (i) of the corollary now follows
from Proposition \ref{UCS}.

Part (ii) of the corollary follows from Theorem \ref{counterthm} and classic 
results of Blum-Hansen \cite{BluHan60}. 
Let $T$ be any strongly mixing invertible transformation 
on $(\X, \calS,\mu )$. Suppose that $H(\C) > 0$ and that $\phi \in \Phi$. 
For $C\in \C$,
\begin{eqnarray*}
\lefteqn{ \int \Big| \frac{1}{n} 
\sum_{i=0}^{n-1} I_C(T_{\phi}^ix) - \mu (C) \Big|^2 \, d\mu } \label{blum1} \\ 
& = & 
\frac{1}{n^2} \sum_{i,j=0}^{n-1} \int 
(I_C(T_{\phi}^i x) - \mu (C)) \, (I_C(T_{\phi}^j x) - \mu (C)) \, d\mu \\[.1in] 
& = & 
\frac{1}{n^2}\sum_{i,j=0}^{n-1} \int 
[I_C(T_{\phi}^i x)I_C(T_{\phi}^j x) - \mu (C)^2] \, d\mu \\[.1in] 
& = &
\frac{1}{n^2}\sum_{i,j=0}^{n-1}
[\mu (C \cap T_{\phi}^{j-i} C ) - \mu (C)^2]. \label{blum2}
\end{eqnarray*}
By Theorem \ref{counterthm} and the Cauchy-Schwartz inequality, 
there exists a dense $G_{\delta}$ subset $G\subset \Phi$ 
such that the first integral does not converge to zero uniformly over $C \in \C$,
and therefore the final sum does not converge to zero uniformly over $C \in \C$. 
From this, a routine argument shows that
\[
\limsup_{n \to \infty} \, \sup_{C \in \C} 
\left| \mu (C \cap T_{\phi}^{-n}C) - \mu (C)^2\right| \, > \, 0
\]
as desired.
$\Box$

\vskip.2in

\subsection{Uniform Weak Mixing}
\label{WkMix}


\vskip.2in

Recall that a measure preserving transformation $T$ on $(\X, \calS, \mu )$ is weak mixing (and hence
ergodic) if for each pair of measurable sets $A, B \in \calS$,
\[
\lim_{n \to \infty} \frac{1}{n} \sum_{i=0}^{n-1} \left| \mu (A \cap T^{-i} B) - \mu(A) \mu (B) \right| \ = \ 0.
\]
The following corollary characterizes uniform weak mixing in terms of entropy.  We omit its proof, which
is substantially similar to that of Corollary \ref{stmix}.

\vskip.2in

\begin{cor}
\label{wkmix} 
Let $T$ be a weak mixing transformation of $(\X, \calS, \mu)$ 
and let $\C \subseteq \calS$.  

\vskip.1in

\begin{enumerate}
\item[(i)] If $H(\C)=0$, then $T$ is uniformly weak mixing on $\C$, in the sense that
\[
\lim_{n \to \infty} \, \sup_{A, B \in \C}
\frac{1}{n} \sum_{i=0}^{n-1} 
\left| \mu(A \cap T^{-i}B) - \mu(A) \mu(B) \right| \, = \, 0.
\]

\vskip.1in

\item[(ii)] Suppose that $T$ is invertible and that $(\X, \calS, \mu)$ is Lebesgue.
If $H(\C) > 0$, then there exists a dense $G_{\delta}$ set $G \subset \Phi$ 
such that for each $\phi \in G$,
\[
\limsup_{n \to \infty} \, \sup_{C \in \C} 
\frac{1}{n} \sum_{i=0}^{n-1} 
\left| \mu (C\cap T_{\phi}^{-i}C) - \mu (C)^2 \right| \, > \, 0.
\]
\end{enumerate}
\end{cor}

\vskip.2in


\vskip.2in


\begin{thebibliography}{1}
\csname bibmessage\endcsname

\bibitem{AdNob12}
\textsc{Adams, T.M.} and \textsc{Nobel, A.B.} (2012)
Uniform approximation of Vapnik-Chervonenkis classes.
\textit{Bernoulli}
\textbf{18:4} 1310-1319.

\bibitem{AdNob10a}
\textsc{Adams, T.M.} and \textsc{Nobel, A.B.} (2010)
Uniform convergence of Vapnik-Chervonenkis classes under 
ergodic sampling.
\textit{Annals of Probability} 
\textbf{38}(4)1345-1367.

\bibitem{AdNob10b}
\textsc{Adams, T.M.} and \textsc{Nobel, A.B.} (2010)
The gap dimension and uniform laws of large numbers for ergodic processes. 
Preprint. arXiv:1007.2964v1

\bibitem{AdSil99}
\textsc{Adams, T.M.} and \textsc{Silva, C.E.} (1999)
$\Z^d$ Staircase actions.
\textit{Ergodic Theory and Dynamical Systems}
\textbf{19} 837-850.

%
%
%
%
%
\bibitem{BluHan60} 
\textsc{Blum, J. R.} and \textsc{Hanson, D. L.}, 
On the mean ergodic theorem for subsequences, 
\textit{Bull. Amer. Math. Soc.}, 
\textbf{66} (1960), 308-311.

%
%
%
%

\bibitem{Dudley99}
\textsc{Dudley, R.M.} (1999) 
\textit{Uniform Central Limit Theorems.}
Cambridge University Press, Cambridge.
MR1720712 (2000k:60050)


\bibitem{FieFri86}
\textsc{Fieldsteel, A.} and \textsc{Friedman, N.A.} (1986)
Restricted orbit changes of ergodic Zd actions to achieve mixing and completely positive entropy.
\textit{Erg. Th. \& Dyn. Sys}
\textbf{6} 505-528.

%

\bibitem{Hal64}
\textsc{Halmos, P.} (1964)
\textit{Measure Theory}, 9th ed., 
D. Van Nostrand Company, New York. MR0033869 (11,504d)

\bibitem{Hal56}
\textsc{Halmos,  P.R.} (1956)
\textit{Lectures on ergodic theory}, Vol. 3,
Amer Mathematical Society.

%
%
%
%
%
%
%
%
%
%
%
%
%
%

\bibitem{Rud78}
\textsc{Rudolph, D.J.} (1978)
The second centralizer of a Bernoulli shift is just its powers.
\textit{Israel Journal of Mathematics}
\textbf{29} (2) 167-178.


%
%
%
%

\bibitem{Tik07}
\textsc{Tikhonov, S. V.} (2007)
A complete metric in the set of mixing transformations. 
\textit{Sb. Math.}
\textbf{198(4):575}


\bibitem{VaaWel96}
\textsc{van der Vaart, A.W.} and \textsc{Wellner, J.A.} (1996)
\textit{Weak Convergence and Empirical Processes}.
Springer-Verlag, New York. 
MR1385671 (97g:60035)

\bibitem{van13}
\textsc{van Handel, R.}(2013)
The universal Glivenko-Cantelli property.
\textit{Probability Theory and Related Fields}
\textbf{155} 911-934.

\bibitem{Vapnik00}
\textsc{Vapnik, V.N.} (2000)
\textit{The nature of statistical learning theory}.  Second edition.
Springer-Verlag, New York. 
MR1719582 (2001c:68110)

%
%
%

\end{thebibliography}
\end{document}